\documentclass{amsart}
\usepackage{amsmath,amsthm,amsfonts,amssymb}

\date{}

\newtheorem{theorem}{Theorem}

\newtheorem{lemma}{Lemma}
\newtheorem{corollary}{Corollary}
\newtheorem{proposition}{Proposition}
\newtheorem{question}{Question}
\newtheorem{problem}{Problem}

\theoremstyle{definition}\newtheorem{remark}{Remark}

\newcommand{\IN}{{\mathbb{N}}}

\newcommand{\IZ}{{\mathbb{Z}}}
\newcommand{\IR}{{\mathbb{R}}}

\newcommand{\IT}{{\mathbb{T}}}

\newcommand{\ind}{\mathrm{ind}}
\newcommand{\Ind}{\mathrm{Ind}}
\newcommand{\w}{\omega}
\newcommand{\e}{\varepsilon}
\newcommand{\cc}{\mathfrak c}

\newcommand{\diam}{\mathrm{diam}}
\newcommand{\concat}{\hat{\phantom{o}}}
\newcommand{\C}{\mathcal C}
\newcommand{\sq}{\mathfrak{sq}}
\newcommand{\dist}{\mathrm{dist}}

\newcommand{\cov}{\mathrm{cov}}

\begin{document}
\title{Packing index of subsets in Polish groups}
\author{Taras Banakh, Nadya Lyaskovska, Du\v{s}an Repov\v{s}}
\thanks{This reseach
was supported the SRA grants P1-029200101-04, J1-9643-0101 and BI-UA/07-08-001} 
\address{Department of Mathematics, Ivan Franko National University of Lviv, Ukraine}
\email{tbanakh@yahoo.com}
\email{lyaskovska@yahoo.com}
\address{Institute of Mathematics, Physics and Mechanics,
and Faculty of Education, University of Ljubljana, P.O.B.2964, Ljubljana, Slovenia}
\email{dusan.repovs@guest.arnes.si}

\begin{abstract} For a subset $A$ of a Polish group $G$, we study the (almost) packing index $\ind_P(A)$ (resp. $\Ind_P(A)$) of $A$, equal to the supremum of cardinalities $|S|$ of subsets $S\subset G$ such that the family of shifts $\{xA\}_{x\in S}$ is (almost) disjoint (in the sense that $|xA\cap yA|<|A|$ for any distinct points $x,y\in S$). Subsets $A\subset G$ with small (almost) packing index are small in a geometric sense. We show that $\ind_P(A)\in \IN\cup\{\aleph_0,\cc\}$ for any $\sigma$-compact subset $A$ of a Polish group. If $A\subset G$ is Borel, then the packing indices $\ind_P(A)$ and $\Ind_P(A)$ cannot take values in the half-interval $[\sq(\Pi^1_1),\cc)$ where $\sq(\Pi^1_1)$ is a certain uncountable cardinal that is smaller than $\cc$ in some models of ZFC. In each non-discrete Polish Abelian group $G$ we construct two closed subsets $A,B\subset G$ with $\ind_P(A)=\ind_P(B)=\cc$ and $\Ind_P(A\cup B)=1$ and then apply this result to show that $G$ contains a nowhere dense Haar null subset $C\subset G$ with $\ind_P(C)=\Ind_P(C)=\kappa$ for any given cardinal number $\kappa\in[4,\cc]$.
\end{abstract}

\maketitle 



\section{Introduction}
Given a Polish group $G$ and a non-empty subset $A\subset G$
with nice descriptive properties, we study 
all
possible values of the packing index
$$
\ind_P(A)=\sup\big\{|S|:\mbox{$S\subset G$ and $\{xA\}_{x\in S}$ is disjoint}\big\}
$$
of $A$, which indicates the smallness of the subset $A$ in a
geometric sense. The papers \cite{BL1}, \cite{BL2}, \cite{L1} and
\cite{L2} are devoted to constructing subsets with a given packing
index. In particular, for every  non-zero cardinal number
$\kappa\le\mathfrak c$ one can easily construct a subset
$A\subset\IR$ with $\ind_P(A)=\kappa$. After discussing those
results on the topological seminar at Wroclaw University the
second author was asked by Krysztof Omiljanowski about possible
restrictions on the packing index $\ind_P(A)$ of subsets $A\subset
\IR$ having good descriptive properties (like being compact
$\sigma$-compact, Borel, measurable or meager). This question was
probably motivated by the well-known fact that the Continuum
Hypothesis (though inresolvable in ZFC) has positive solution in
the realm of Borel subsets of the real line: each uncountable
Borel subset $A\subset\IR$ has cardinality $\mathfrak c$ of
continuum.

In this paper we shall give several partial answers to 
Omiljanowski's question.
On the
one hand, we show that $\sigma$-compact subsets $A$ in Polish groups cannot have
an intermediate packing index $\aleph_0<\ind_P(A)<\mathfrak c$. For a Borel subset $A$
of a Polish group we have a weaker result: $\ind_P(A)$ cannot take the value in the
interval $\sq(\Pi^1_1)\le\ind_P(A)<\cc$ where $\sq(\Pi^1_1)$ stands for
the smallest cardinality $\kappa$ such that each coanalytic subset
$X\subset 2^\w\times 2^\w$ contains a square $S\times S$ of size
$|S\times S|=\cc$ provided $X$ contains a square of size $\ge\kappa$.
The value of the small uncountable cardinal $\sq(\Pi^1_1)$ is not completely
determined by ZFC Axioms: both the equality $\sq(\Pi^1_1)=\cc$ and the strict
inequality $\sq(\Pi^1_1)<\cc$ are consistent with Martin Axiom, see \cite{Sh}.

On the other hand, for every infinite cardinal number $\kappa\le\mathfrak c$ in
each non-discrete Polish Abelian group $G$ we shall construct a nowhere dense Haar
null subset $A\subset G$ with $\ind_P(A)=\Ind_{P}(A)=\kappa$. Here
$$\Ind_{P}(A)=\sup\{|S|:\mbox{$S\subset G$ and $\{xA\}_{x\in S}$ is almost disjoint}\}$$
is the {\em almost packing index} of $A$. In the above definition,
a family of shifts $\{xA\}_{x\in S}$ is defined to be {\em almost
disjoint} if $|xA\cap yA|<|A|$ for all distinct $x,y\in S$. 

To construct the nowhere dense Haar null subset $A\subset G$ with a given (almost) packing index, in each non-discrete Polish
Abelian group $G$ we first construct a closed nowhere dense Haar null
subset $C\subset G$ with $\Ind_{P}(C)=1$. The set $C$, being
nowhere dense and Haar null, is small in the sense of category and
measure, but is large in the geometric sense because for any two
distinct points $x,y\in G$ the shifts $xC$ and $yC$ have
intersection of cardinality $|xC\cap yC|=|G|$. In particular,
$CC^{-1}=G$ and thus $C$ is a closed nowhere dense Haar null
subset that algebraically generates the group $G$. This result can
be seen as an extension of a result of S.Solecki \cite{So} who
proved that each non-locally compact Polish Abelian group $G$ is algebraically
generated by a nowhere dense subset.  Also it extends some results
of \cite[\S13]{BP}. In fact, the closed Haar null subset $C\subset
G$ with $\Ind_P(C)=1$ is constructed as the union $C=A\cup B$ of
two closed subsets $A,B\subset G$ with $\ind_P(A)=\ind_P(B)=\cc$.
This shows that the packing index is highly non-additive.
\medskip

{\bf Notation.} By $\w$ we denote the first infinite ordinal,
 $\IN=\w\setminus\{0\}$ stands for the set of positive integers.
 Cardinals are identified with the initial ordinals of given cardinality;
 $\mathfrak c$ stands for the cardinality of continuum. All topological groups $G$
 considered in this paper will be supplied with an invariant metric $\rho$ generating
 the topology of $G$. By $e$ we denote the identity element of $G$. For a point $x\in G$
 and a real number $r$ by $B(x,r)=\{g\in G:\rho(g,x)<r\}$ we denote the open $r$-ball
  centered at $x$. Also for $x\in G$ we put $\|x\|=\rho(x,e)$. The invariantness of $\rho$
  implies $\|x^{-1}\|=\|x\|$ and $\|xy\|\le\|x\|+\|y\|$ for all $x,y\in G$.

\section{The packing indices of $\sigma$-compact sets in Polish groups}

In this section we show that the packing index of a $\sigma$-compact subset in a Polish groups cannot take an intermediate value between $\w$ and $\cc$.

First we prove a useful

\begin{lemma}\label{aaa} Let $A$ be a subset of a
Polish group $G$. If\/ $\ind_P(A)<\cc$, then the closure of $AA^{-1}$
 contains a neighborhood of the neutral element $e$ of $G$.
\end{lemma}

\begin{proof}
Fix any complete metric $\rho$ generating the topology of the Polish group $G$.
Assuming that $\overline{AA^{-1}}$ is not a neighborhood of $e$,  we shall construct
a perfect subset $K\subset G$ such that $(xA)_{x\in K}$ is disjoint, which will imply
that $\ind_P(A)=|K|=\cc$.

Taking into account that the closed subset $\overline{AA^{-1}}$ is not a neighborhood
of $e$ in $G$, for any open neighborhood
$U$ of $e$ we can find a point $b\in U\setminus \overline{AA^{-1}}$ and an open
neighborhood $V$ of
$e$ such that $V^{-1}bV\subset U\backslash AA^{-1}$.

Using this fact, by induction construct  a sequence $(b_n)_{n\in\w}$ of points in
$G$ and sequences $(U_n)_{n\in\w}$ and $(V_n)_{n\in\w}$ of open neighborhoods of $e$
in $G$ such that
\begin{enumerate}
\item $b_n\in U_n=U_n^{-1}$;
\item $V_{n+1}^{-1}b_nV_{n+1}\cap AA^{-1}=\emptyset$;
\item $b_n\notin V_{n+1}V_{n+1}^{-1}$;
\item $\diam_\rho(bV_{n+1})<2^{-n}$ for any
point\newline $b\in B_n=\big\{b_0^{\e_0}\cdots b_n^{\e_n}:\e_0,\dots,\e_n\in\{0,1\}\big\}$;
\item $U_{n+1}^3\subset V_{n+1}\subset U_n$.
\end{enumerate}

Define a map $f:\{0,1\}^\w\to G$ assigning to each infinite binary
sequence $\vec\e=(\e_i)_{i\in\w}$ the infinite product
$$f(\vec \e)=\prod_{i=0}^\infty b_i^{\e_i}=\lim_{n\to\infty} f_n(\vec \e)$$
 where $f_n(\vec \e)=\prod_{i=0}^n b_i^{\e_i}$. Let us show that the latter limit exists.
 It suffices to check that $(f_n(\vec\e))_{n\in\w}$ is a Cauchy sequence in $(G,\rho)$.

The condition (5) implies that $U_{i+1}^2\subset U_i$ for all $i$.
This can be used as the inductive step in the proof of the inclusion
$U_n\cdots U_m\subset U_n^2$ for all $m\ge n$. Then for every $m\ge n$
$$f_m(\vec \e)\in f_n(\vec \e)\cdot U_{n+1}\cdots U_m\subset f_n(\vec\e)\cdot U^2_{n+1}
\subset f_n(\vec\e)V_{n+1}$$
and thus
$$\rho(f_m(\vec\e),f_n(\vec
\e))\le\diam_\rho(f_n(\vec\e)\cdot V_{n+1})<2^{-n}$$
 by the
condition (4). Therefore, the sequence $(f_n(\vec\e))_{n\in\w}$ is
Cauchy and the limit $f(\vec\e)=\lim_{n\to\infty}f_n(\vec\e)$
exists. Moreover, the upper bound $\rho(f_m(\vec
\e),f_n(\vec\e))\le 2^{-n}$ implies that the map $f:\{0,1\}^\w\to
G$ is continuous. On the other hand, the inclusions
$f_m(\vec\e)\in f_n(\vec\e)\cdot U_n^2$, $m\ge n$, imply that
$$f(\vec\e)\in f_n(\vec\e)\cdot \overline{U_n^2}\subset
f_n(\vec\e)\cdot U_n^3\subset f_n(\vec\e)\cdot V_{n+1}.$$

This inclusion will be used for the proof of the injectivity of
$f$. We shall prove a little bit more: for any distinct vectors
$\vec\e$ and $\vec \delta$ in $\{0,1\}^\w$ we get $f(\vec \e)A\cap
f(\vec\delta)A=\emptyset$. Find the smallest number $n\in\w$ such
that $\e_n\ne\delta_n$. We lose no generality assuming that
$\delta_n=0$ and $\e_n=1$. It follows that $f(\vec\e)\in f_n(\vec
\e)U_{n+1}^3=f_{n-1}(\vec\e)b_nV_{n+1}$ while
$f(\vec\delta)=f_n(\vec\delta) V_{n+1}=f_{n-1}(\delta) \cdot
e\cdot V_{n+1}=f_{n-1}(\vec\e)V_{n+1}$. Then
$$\big(f(\vec\delta)\big)^{-1}f(\vec\e)\in V_{n+1}^{-1}b_n
V_{n+1}\subset G\setminus AA^{-1}$$ by the condition (2) and hence
$f(\vec\e)A\cap f(\vec\delta)A=\emptyset$.

Thus the family $(xA)_{x\in K}$ is disjoint where $K=f(\{0,1\}^\w)$.
The injectivity of $f$ implies that $\ind_P(A)\ge|K|=\cc$.
\end{proof}

Now we can prove the main result of this section.

\begin{theorem}\label{t1} If $A$ is a $\sigma$-compact
subset $A$ in a Polish group $G$, then
$\ind_P(A)\in \IN\cup\{\aleph_0,\cc\}$. Moreover, if the set $A$ is compact, then
\begin{enumerate}
\item $\ind_P(A)=\cc$ if $G$ is not locally compact;
\item $\ind_P(A)\in \{\aleph_0,\cc\}$ if $G$ is locally compact but not compact;
\item $\ind_P(A)\in \IN\cup\{\cc\}$ if $G$ is compact.
\end{enumerate}
\end{theorem}

\begin{proof} If $A$ is $\sigma$-compact, then so is the set
$AA^{-1}=\{xy^{-1}:x,y\in A\}$ and then the set $(G\setminus AA^{-1})\cup\{e\}$
is a $G_\delta$-set in $G$. In its turn,
the subset
$$X=\{(x,y)\in G\times G:y^{-1}x\in (G\setminus AA^{-1})\cup\{e\}\}$$
is of type $G_\delta$ in $G\times G$, being the preimage of the
$G_\delta$-subset $(G\setminus AA^{-1})\cup\{e\}$ under the
continuous map $g:G\times G\to G$, $g:(x,y)\mapsto y^{-1}x$.

Assuming that $\ind_P(A)>\aleph_0$, we could find an uncountable
subset $S\subset G$ with disjoint family $\{xA\}_{x\in S}$, which
implies that $S\times S\subset X$. Since the Polish space $X$
contains the uncountable square $S\times S$, we can apply the
Shelah's result \cite[1.14]{Sh}  to conclude that $X$ contains the
square $P\times P$ of a perfect subset $P\subset G$ (the latter
means that $P$ is closed in $G$ and has no isolated point). It
follows from $P\times P\subset X$ that the family $\{xA\}_{x\in
P}$ is disjoint and thus $\cc=|P|\le\ind_P(A)\le|G|=\cc$.
\smallskip

Now assuming  that $A\subset G$ is compact we shall prove the items (1)--(3)
of the theorem. The compactness of $A$ implies the compactness of $AA^{-1}$.
If $AA^{-1}$ is not a neighborhood of $e$, then we can apply Lemma~\ref{aaa} to
conclude that $\ind_P(A)=\cc$. This is so if the group $G$ is not locally compact.
 So, next we assume that $AA^{-1}$ is a neighborhood of $e$. In this case the group
 is locally compact and we can take a neighborhood $U\subset G$ of $e$ with
 $UU^{-1}\subset AA^{-1}$.

 Then for every $B\subset G$ with $B^{-1}B\cap AA^{-1}=\{e\}$ we get
 $B^{-1}B\cap UU^{-1}=\{e\}$, which implies that the family $(xU)_{x\in B}$
 is disjoint and the set $B$ is at most countable, being discrete in the Polish
 space $G$. This proves the upper bounds $\ind_P(A)\le\aleph_0$.

If the group $G$ is not compact, then using the Zorn Lemma, we can find a maximal
set $B\subset G$ with $B^{-1}B\cap AA^{-1}=\{e\}$. We claim that $BAA^{-1}=G$. Assuming
the converse, we could find a point $b\in G\setminus BAA^{-1}$. Then the set $bA$ is disjoint
from the set $BA$ and hence we can enlarge the set $B$ to the set $\tilde B=B\cup\{b\}$ such
 that $(xA)_{x\in\tilde B}$ is disjoint. The latter is equivalent to
 $\tilde B^{-1}\tilde B\cap AA^{-1}=\{e\}$ and this contradicts the maximality of $B$.
The compactness of $AA^{-1}$ and non-compactness of $G=BAA^{-1}$ implies that $B$
is infinite and thus $\ind_P(A)\ge|B|\ge\aleph_0$. This completes the proof of the
second item of the theorem.

To prove the third item, assume that $G$ is compact. In this case $G$ carries a Haar
measure $\mu$ which is a unique probability invariant $\sigma$-additive Borel measure
on $G$. If $AA^{-1}$ is not a neighborhood of $e$, then  $\ind_P(A)=\cc$ by a preceding
case. So we assume that $AA^{-1}$ is a neighborhood of $e$ and take another neighborhood
$U$ of $e$ with $UU^{-1}\subset AA^{-1}$. Since finitely many shifts of $U$ cover the group
$G$, we get $\mu(U)>0$. Now given any subset $B\subset G$ with $B^{-1}B\cap AA^{-1}=\{e\}$,
we get $B^{-1}B\cap UU^{-1}=\{e\}$. The  latter equality implies that the family $(xU)_{x\in B}$
is disjoint and then $1=\mu(G)\ge\mu(BU)=|B|\mu(U)$ implies that $|B|\le1/\mu(U)$.
Consequently, the packing index $\ind_P(A)\le1/\mu(U)$ is finite.
\end{proof}

  In light of this theorem two open questions arise naturally.

\begin{question}\label{q1} Is there a compact group $G$ and a $\sigma$-compact
subset $A$ with $\ind_P(A)=\aleph_0$?
\end{question}

\begin{question}\label{q2} Is there a Polish group $G$ and a Borel
subset $A\subset G$ with $\aleph_0<\ind_P(A)<\cc$?
\end{question}

The latter question does not reduce to the $\sigma$-compact case because of the following
example (in which $\IT=\IR/\IZ$ stands for the circle).

\begin{proposition} The countable product $G=\IT^\w$ contains a $G_\delta$-subset
$A\subset G$ such that
$\ind_P(A)=\cc$ but each $\sigma$-compact subset $B\subset \IT^\w$ containing $A$
has $\ind_P(B)<\aleph_0$.
\end{proposition}

\begin{proof} Let $q:\IZ\to\IT$ denote the quotient map, $J=q\big([0,\frac12)\big)$
be the half-circle,  $I=\overline{J}=q([0,\frac12])$ be its closure, and $D=q(\{0,\frac12\})$
be two opposite points on $\IT$. It is clear that $D^{-1}D\cap JJ^{-1}=\{e\}$ while
$I\cdot I^{-1}=\IT$.

It follows that $A=J^\w$ is a $G_\delta$-subset of $\IT^\w$ with $\ind_P(A)=|D^\w|=\cc$
because $(D^\w)^{-1}D^\w\cap AA^{-1}=\{e\}$.

Now given any $\sigma$-compact subset $B\supset A$ in $\IT^\w$, we should check that
$\ind_P(B)<\aleph_0$. Replacing $B$ by $B\cap I^\w$, if necessary, we can assume that
 $B\subset I^\w$. Since $B\subset I^\w$ contains the dense $G_\delta$-subset $J^\w$ of $I^\w$,
  the standard application of the Baire Theorem yields an non-empty open subset
  $U\subset I^\w$ with $U\subset B$.
We lose no generality assuming that $U$ is of basic form $U=V\times I^{\w\setminus n}$
for some $n\in\w$ and some open set $V\subset I^n$. Observe that
$$U^{-1}U=VV^{-1}\times I^{\w\setminus n}(I^{\w\setminus n})^{-1}=VV^{-1}\times
\IT^{\w\setminus n}$$
is an open neighborhood of $e$ in $\IT^\w$. Consequently,
$BB^{-1}\supset UU^{-1}$ is also an open neighborhood of $e$ in
$\IT^\w$. Proceeding as in the proof of the first item of
Theorem~\ref{t1}, we can see that $$\ind_P(B)\le
1/\mu(V\times\IT^{\w-n})<\aleph_0.$$
\end{proof}

\section{The packing indices of analytic sets in Polish groups}

In this section we shall give a partial answer to
Question~\ref{q2} related to the packing indices of Borel subsets
in Polish groups. It is well-known that each Borel subset of a
Polish space is analytic. We recall that a metrizable space $X$ is
{\em analytic} if $X$ is a continuous image of a Polish space. A
space $X$ is {\em coanalytic }if for some Polish space $Y$
containing $X$ the complement $Y\setminus X$ is analytic. The
classes of analytic and coanalytic spaces are denoted by
$\Sigma^1_1$ and $\Pi^1_1$, respectively. It is well-known that
the intersection $\Delta^1_0=\Sigma^1_1\cap \Pi^1_1$ coincides
with the class of all absolute Borel (metrizable separable)
spaces. By $\Sigma^0_2$ and $\Pi^0_2$ we denote the classes of
$\sigma$-compact and Polish spaces, respectively.

We shall say that a subset $X\subset 2^\w\times2^\w$ contains a square of size
 $\kappa$ if there is a subset $A\subset 2^\w$ with $A\times A\subset X$ and
 $|A\times A|=\kappa$. Given a class $\C$ of spaces denote by $\sq(\C)$ the smallest
 cardinal $\kappa$ such that each subspace $X\in\C$ of $2^\w\times 2^\w$ that contains
 a square of size $\kappa$ contains a square of size $\cc$. The Shelah's
 result \cite{Sh} (applied in the proof of Theorem~\ref{t1}) guarantees
 that $\sq(\Pi^0_2)=\aleph_1$. For other descriptive classes $\C$ the value
 $\sq(\C)$ is not so definite and depends on Set-Theoretic Axioms.
In particular, the Continuum Hypothesis implies that
$\sq(\Sigma^0_2)=\sq(\Sigma^1_1)=\sq(\Pi^1_1)=\cc$. On one hand,
the strict inequality $\sq(\Pi^1_1)<\cc$ is consistent with ZFC+MA, see
\cite[1.9, 1.10]{Sh}. However, there is a substantial difference between the
classes $\Pi^0_2$ and $\Sigma^0_2$ of Polish and $\sigma$-compact
spaces. By \cite{Sh} each Polish space $X\subset 2^\w\times 2^\w$
containing an uncountable square contains a Perfect square. On the
other hand, there is a ZFC-example of a $\sigma$-compact subspace
$X\subset 2^\w\times 2^\w$ that contains a square of size
$\aleph_1$ but not the perfect one, see \cite{K}.

\begin{proposition} Let $A$ be an analytic subset of a Polish group $G$.
If $\ind_P(A)\ge\sq(\Pi^1_1)$, then $\ind_P(A)=\cc$.
\end{proposition}

\begin{proof} Using the fact that each Polish space is a continuous one-to-one image
of a zero-dimensional Polish space, we can show that $\sq(\Pi^1_1)$ coincides with the
smallest cardinal $\kappa$ such that for any Polish space $X$ a coanalytic subset
$C\subset X\times X$ contains a square of size $\cc$ provided $C$ contains a square of
size $\ge \kappa$.

Given an analytic subset $A$ of a Polish group $G$ we can see that
both the sets $AA^{-1}$ and $AA^{-1}\setminus\{e\}$ are analytic
and thus the set $C=\{(x,y)\in G\times G:y^{-1}x\notin
AA^{-1}\setminus\{e\}\}$ is coanalytic.

Assuming that $\ind_P(A)\ge\sq(\Pi^1_1)$, we could find a subset $S\subset G$
of size $|S|\ge\sq(\Pi^1_1)$ such that the family $\{xA\}_{x\in S}$ is disjoint.
The latter is equivalent to $S^{-1}\cdot S\subset G\setminus (AA^{-1}\setminus\{e\})$
and to $S\times S\subset C$. By the equivalent definition of $\sq(\Pi^1_1)$ (with $2^\w$
replaced by any Polish space), the coanalytic subset $C\subset G\times G$ contains a
square $K\times K$ of size $\cc$ (because it contains the square $S\times S$ of cardinality
$|S\times S|\ge\sq(\Pi^1_1)$).
 It follows from $K\times K\subset C$ that the family $\{xA\}_{x\in K}$ is disjoint
 and thus $\ind_P(A)\ge|K|=\cc$.
\end{proof}

A similar result holds for the almost packing index
$$\Ind_P(A)=\sup\{|S|:\mbox{$S\subset G$ and $\{xA\}_{x\in S}$ is almost disjoint$\}$ of $A$}.$$
 We recall that $\{xA\}_{x\in S}$ is {\em almost disjoint} if $|xA\cap yA|<|A|$ for any
 distinct points $x,y\in S$.

In the proof of the following theorem we shall use a known fact of the Descriptive Set
Theory saying that for a Borel subset $A\subset X\times Y$ in the product of two Polish spaces
the set $\{y\in Y:|A\cap (X\times\{y\})|\le\aleph_0\}$ is coanalytic in $Y$, see \cite[18.9]{Ke}.

\begin{proposition} Let $A$ be a Borel subset of a Polish group $G$.
If\/ $\Ind_P(A)\ge\sq(\Pi^1_1)$, then $\Ind_P(A)=\cc$.
\end{proposition}

\begin{proof} It follows from $\Ind_P(A)\ge\sq(\Pi^1_1)>\aleph_0$ that the space $G$
is uncountable. If $A$ is countable, then trivially, $\Ind_P(A)=\ind_P(A)=\cc$.

So we assume that $A$ is uncountable.
First we show that the subset $C=\{x\in G:|A\cap xA|\le\aleph_0\}$ is coanalytic.
Consider the homeomorphism $h:G\times G\to G\times G$, $h:(x,y)\mapsto (x,y^{-1}x)$, and
the Borel subset $B=h(A\times A)\subset G\times G$. Since
$C=\{z\in G:|B\cap (G\times\{z\})|\le\aleph_0\}$, we may apply the mentioned result
\cite[18.9]{Ke} to conclude that the set $C$ is coanalytic. Then the set
$D=\{(x,y)\in G\times G:y^{-1}x\in C\}$ is coanalytic as the preimage of the coanalytic
subset under a continuous map between Polish spaces.

Assuming that $\Ind_P(A)\ge\sq(\Pi^1_1)$, we could find a subset $S\subset X$ such
that $|S|\ge\sq(\Pi^1_1)$ such that the family $\{xA\}_{x\in S}$ is almost \mbox{disjoint}.
Then for any distinct $x,y\in S$ the intersection $xA\cap yA$, being a Borel subset of
cardinality $|xA\cap yA|<|A|\le\cc$, is at most countable. Consequently, $y^{-1}x\in C$ and
 thus $S\times S\subset D$. By the equivalent definition of $\sq(\Pi^1_1)$, the coanalytic set
  $D$ contains a square $K\times K$ of size $\cc$. It follows from $K^{-1}K\subset C$ that the
  family $\{xA\}_{x\in K}$ is almost disjoint. Consequently, $\cc=|K|\le\Ind_P(A)\le |G|=\cc$.
\end{proof}

\begin{question} Let $A$ be a Borel subset of a Polish group $G$. Is there an at most
countable subset $C\subset A$ such that  $\ind_P(A\setminus C)=\Ind_P(A)$?
\end{question}

The other problem concerns the cardinals $\sq(\C)$ for various descriptive classes
 $\C$. If such a class $\C$ contains the square of a countable metrizable space,
  then $\aleph_1\le\sq(\C)\le\cc$ and thus $\sq(\C)$ falls into the category of the
  so-called small uncountable cardinals, see \cite{V}. However unlike to other typical
   small uncountable cardinals, $\sq(\C)$ does no collapse to $\cc$ under the Martin Axiom, see \cite{Sh}.

\begin{problem} Explore possible values and inequalities between
classical small uncountable cardinals and the cardinals $\sq(\C)$ for various
descriptive classes $\C$.
\end{problem}

\section{Relation of the packing index to other notions of smallness}

Taking into account that a subset $A$ with large packing index $\ind_P(A)$ is
geometrically small, it is natural to consider the relation of the packing index to some
other known concepts of smallness, in particular, the smallness in the sense of Baire category
and the measure.

We recall that a subset $A$ of a topological space $X$ is {\em meager} if $A$
can be written as the countable union of nowhere dense subsets. We shall need
the following classical fact.

\begin{proposition}[Banach-Kuratowski-Pettis]\label{bkp} For any analytic non-meager
subset $A$ of a Polish group $G$ the set $AA^{-1}$ contains a neighborhood of the neutral
element of $G$.
\end{proposition}

A similar result holds for analytic subsets that are not Haar
null. We recall that a subset $A$ of a topological group $G$ is
called {\em Haar null} if there is a Borel probability measure
$\mu$ on $G$ such that $\mu(xAy)=0$ for all $x,y\in G$. This
notion was introduced by J.Christensen \cite{C} and thoroughly
studied  in \cite{THJ}. In particular, a subset $A$ of a
locally compact group $G$ is Haar null if and only if $A$ has zero
Haar measure. Yet, Haar null sets exists in non-locally compact
groups (admitting no invariant measure).

\begin{proposition}[Chistensen]\label{ch} If an analytic subset $A$ of a
Polish group $G$ is not Haar null, then $AA^{-1}$ contains a neighborhood of
the neutral element of $G$.
\end{proposition}

We shall use those propositions to prove

\begin{theorem}\label{megHar} Let $A$ be an analytic subset of a Polish group $G$.
If $\ind_P(A)>\aleph_0$, then $A$ is meager and Haar null.
\end{theorem}

\begin{proof} Otherwise, we can apply Propositions~\ref{bkp} or \ref{ch}
to conclude that $AA^{-1}$ contains a neighborhood $U$ of the neutral element $e$ of $G$.
 Find another neighborhood $V\subset G$ of $e$ with $VV^{-1}\subset U\subset AA^{-1}$.

Since $\ind_P(A)>\aleph_0$ there is an uncountable subset $S\subset X$ such that the
 family $\{xA\}_{x\in S}$ is disjoint, which is equivalent to $S^{-1}S\cap AA^{-1}=\{e\}$.
 It follows from the choice of $V$ that $S^{-1}S\cap VV^{-1}\subset S^{-1}S\cap AA^{-1}=\{e\}$
  and thus the family $\{xV\}_{x\in S}$ is disjoint. Since $V$ is an open neighborhood of $e$,
  the set $S$, being discrete in $G$, is at most countable. This contradiction completes the
  proof.
\end{proof}

\section{The packing index and unions}

It is known that the countable union of meager (resp. Haar null) subsets of a Polish group
is meager (resp. Haar null). In contrast, the union of two subsets $A,B\subset G$ with
large packing index need not have large packing index. A simplest example is given by the
sets $A=\IR\times\{0\}$ and $B=\{0\}\times\IR$ on the plane $\IR^2$. They have packing indices
 $\ind_P(A)=\ind_P(B)=\cc$ but $\ind_P(A\cup B)=1$. In fact, this situation is typical.
 According to \cite{L2}, each infinite group $G$ contains two sets $A,B\subset G$ such that
 $\ind_P(A)=\ind_P(B)=|G|$ but $\Ind_P(A\cup B)=1$. The following theorem is a topological
 version of this result.

\begin{theorem}\label{t3} Each non-discrete Polish Abelian group $G$ contains
two closed subsets $A,B\subset G$ such that $\ind_P(A)=\ind_P(B)=\cc$ and
$\Ind_P(A\cup B)=1$.
\end{theorem}

\begin{proof} Fix an invariant metric $\rho$ generating the topology of $G$.
This metric is complete because $G$ is Polish.  Since $G$ is Abelian, we use the
 additive notation for the group operation on $G$. The neutral element of $G$ will be
 denoted by $0$.

We define a subset $D$ of $G$ to be {\em $\varepsilon$-separated}
if $\rho(x,y)\ge\varepsilon$ for any distinct points $x,y\in D$. By the Zorn lemma,
each $\e$-separated subset can be enlarged to a maximal $\e$-separated subset of $G$.

Put $\e_{-1}=\e_0=1$ and choose a maximal $2\e_0$-separated subset $H_0\subset G$
containing zero.
Proceeding by induction we shall define a sequence $(h_n)_{n\in\IN}\subset G$ of points,
a sequence $(\varepsilon_n)_{n\in \omega}$ of positive
real numbers  and a sequence $(H_n)_{n\in \omega}$ of subsets of $G$ such that
for every $n>0$
\begin{enumerate}
\item [(i)]$B(0,\e_{n-1})\setminus B(0,33\e_{n})$ is not empty and $\e_n<2^{-6}\e_{n-1}$;
\item [(ii)] $\|h_n\|=5\varepsilon_{n}$,
\item [(iii)] $H_n\supset\{0,h_n\}$ is a maximal $2\varepsilon_{n}$-separated subset
in $B(0, 8\varepsilon_{n-1})$.
\end{enumerate}

It follows from (i) that the series $\sum_{n\in\w}\e_n$ is convergent and thus for
any sequence $x_n\in H_n$, $n\in\w$, the series $\sum_{n\in\w}x_n$ is convergent
(because $\|x_n\|<8\e_{n-1}$ for all $n\in\IN$ according to (iii)). So it is legal
to consider the sets of sums
$$
\begin{aligned}
\Sigma_0&=\big\{\sum_{n\in \omega}x_{2n}:(x_{2n})_{n\in\w}\in\prod_{n\in\w} H_{2n}\big\},\\
\Sigma_1&=\big\{-\sum_{n\in \omega}x_{2n+1}:(x_{2n+1})_{n\in\w}\in\prod_{n\in\w} H_{2n+1}\big\}.
\end{aligned}$$

Let $A$ and $B$ be the closures of the sets $\Sigma_0$ and $\Sigma_1$ in $G$. It remains
 to prove that the sets $A,B$ have the desired properties: $\ind_P(A)=\ind_P(B)=\cc$ and
 $\Ind_P(A\cup B)=1$. This will be done in three steps.
\smallskip

1. First we prove that $\ind_P(A)=\cc$. By Lemma~\ref{aaa}, this equality will follow
as soon as we check  that $\overline{A-A}$ is not a neighborhood of the neutral element
$0$ in $G$.
It suffices for every $k\in\w$ to find a point a point
 $g\in B(0,\e_{2k})\setminus \overline{AA^{-1}}$. By condition (i), there is a point
 $g\in G$ with $33\e_{2k+1}\le \|g\|<\e_{2k}$. We claim that
 $g\notin \overline{A-A}=\overline{\Sigma_0-\Sigma_0}$. More precisely,
 $$\dist(g,\overline{A-A})=\dist(g,\Sigma_0-\Sigma_0)\ge\min\{\e_{2k+1},\e_{2k}/2\}.$$

Take any two distinct points $x,y\in \Sigma_0$ and find infinite sequences
$(x_{2n})_{n\in\w},(y_{2n})_{n\in\w}\in\prod_{n\in\w}H_{2n}$ with $x=\sum_{n\in\w}x_{2n}$
 and $y=\sum_{n\in\w}y_{2n}$.

Let $m=\min\{n\in\w:x_{2n}\ne y_{2n}\}$. If $m\ge k+1$, then
$$
\begin{aligned}
\|x-y\|=&\,\|\sum_{n\ge m}x_{2n}-y_{2n}\|\le\sum_{n\ge m}\|x_{2n}\|+\|y_{2n}\|\le\\
\le&\,2\sum_{n\ge m}8\e_{2n-1}\le 32\,\e_{2m-1}\le32\,\e_{2k+1}<\|g\|-\e_{2k+1}
\end{aligned}$$ and hence $\rho(x-y,g)\ge\e_{2k+1}$.

If $m\le k$, then
$$
\begin{aligned}
\|x-y\|=&\,\|(x_{2m}-y_{2m})+\sum_{n>m}(x_{2n}-y_{2n})\|\ge\\
\ge&\, \|x_{2m}-y_{2m}\|-\sum_{n>m}(\|x_{2n}\|+\|y_{2n}\|)\ge\\
\ge&\, 2\e_{2m}-32\e_{2m+1}\ge \frac32\e_{2m}>\|g\|+\frac12\e_{2m}
\end{aligned}$$ and again $\rho(x-y,g)\ge\frac12\e_{2m}\ge\frac12\e_{2k}$.
\smallskip

2. In the same manner we can prove that $\ind_P(B)=\cc$.
\smallskip

3. It remains to prove that $\Ind_P(A\cup B)=1$. First we recall some standard notation.
Denote by $2^{<\w}=\bigcup_{n\in\w}2^n$ the set of finite binary sequences. For any sequence
$s=(s_0,\dots,s_{n-1})\in 2^{<\w}$ and $i\in2=\{0,1\}$ by $|s|=n$ we denote the length
of $s$ and by
$s\hat{\phantom{o}}i=(s_0,\dots,s_{n-1},i)$ the concatenation of $s$ and $i$. For a finite
or infinite binary sequence $s=(s_i)_{i<n}$ and $l\le n$ let $s|l=(s_i)_{i<l}$.
The set $2^{\w}$ is a tree with respect to the partial order: $s\le t$ iff  $s=t|l$ where
$l=|s|\le|t|$.

The equality $\ind_P(A\cup B)=1$ will follow as soon as we prove that $|A\cap (g+B)|\ge\cc$
for all $g\in G$.
We shall construct a sequence of points $\{x_s\}_{s\in 2^{<\omega}}$ such that
for every sequence $s\in 2^{<\w}$ the following conditions hold:
\begin{enumerate}
\item $x_s\in H_{|s|}\subset B(0,8\e_{|s|-1})$;
\item $\|x_{s\concat0}-x_{s\concat1}\|>\e_n$;
\item $\|g-\sum_{t\le s}x_t)\|<7 \e_{|s|}$.
\end{enumerate}

We start choosing a point $x_\emptyset\in H_0$  with $\rho(x_\emptyset,g)<2\e_{-1}=2\e_0$.
Such a point $x_\emptyset$ exists because $H_0$ is a maximal $(0,2\e_0)$-separated set in $G$.
 Next we proceed by induction. Suppose that for some $n$ the points $x_s$, $s\in 2^{<n}$,
 have been constructed. Given a sequence $s\in 2^{n-1}$ we need to define the points
 $x_{s\concat0}$ and $x_{s\concat1}\in H_n$. Let $g_s=g-\sum_{t\le s}x_t$.
 Since
$\|g_s+h_n\|\le\|g_s\|+\|h_n\|<7\e_{n-1}+5\e_n<8\e_{n-1}$ and  $H_n$ is a maximal
$2\e_n$-separated subset in $B(0,8\e_{n-1})$, there are two points
$x_{s\concat0},x_{s\concat1}\in H_n$ with $\rho(g_s,x_{s\concat0})<2\e_n$ and
$\rho(g_s+h_n,x_{s\concat1})<2\e_n$. The condition (2) follows from
$$
\begin{aligned}
\|x_{s\concat0}-x_{s\concat1}\|&\ge \|g_s-(g_s+h_n)\|-\|g_s-x_{s\concat0}\|-\|g_s+h_n-x_{s\concat1}\|>\\
&>5\e_n-2\e_n-2\e_n=\e_n.
\end{aligned}$$ The condition (3) follows from the estimates $$
\|g-\sum_{t\le s\concat0}x_t\|=\|g-x_{s\concat0}-\sum_{t\le s}x_t\|=\|g_s-x_{s\concat0}\|<2\e_n=2\e_{|s\concat0|}$$and
$$
\begin{aligned}
\|g-\sum_{t\le s\concat1}x_t\|=&\,\|g-x_{s\concat1}-\sum_{t\le s}x_t\|=\|g_s+h_n-x_{s\concat1}-h_n\|\le\\
\le&\,\|g_s+h_n-x_{s\concat1}\|+\|h_n\|<2\e_n+5\e_n=7\e_{|s\concat1|}.
\end{aligned}$$

After completing the inductive construction, we can use the condition (3) to see that for
every infinite binary sequence $s\in 2^\w$ we get

$$g=\sum_{n\in\w}x_{s|n}=\sum_{n\in\w}x_{s|2n}+\sum_{n\in\w}x_{s|2n+1}.$$
We claim that the set $$D_0=\big\{\sum_{n\in\w}x_{s|2n}:s\in
2^\w\big\}$$ lies in the intersection $A\cap(g+B)$. It is clear
that $D_0\subset\Sigma_0\subset A$. To see that $D_0\subset g+B$,
take any point $x\in D_0$ and find an infinite binary sequence
$s\in 2^\w$  with $x=\sum_{n\in\w}x_{s|2n}$. Then
$$x=\sum_{n\in\w}x_{s|2n}+\sum_{n\in\w}x_{s|2n+1}-\sum_{n\in\w}x_{s|2n+1}\in
g+\Sigma_1\subset g+B.$$

 It remains to prove that $|D_0|\ge \cc$. Note that the set $D_0$, being a continuous
 image of the Cantor cube $2^\w$, is compact. Now the equality $|D_0|=\cc$ will follow
 as soon as we check that $D_0$ has no isolated points. Given any sequence $s\in 2^\w$
 and $\delta>0$ we should find a sequence $t\in2^\w$ such that

 $$0<\|\sum_{n\in\w}x_{s|2n}-\sum_{n\in\w}x_{t|2n}\|<\delta.$$
Find even number $2m\in\w$ such that $\sum_{n\ge m}20\e_{2n-1}<\delta$ and take any
sequence $t\in2^\w$ such that $t|2m-1=s|2m-1$ but $t|2m\ne s|2m$.  Then $$
\begin{aligned}
\|\sum_{n\in\w}x_{s|2n}-\sum_{n\in\w}x_{t|2n}\|=&\,\|\sum_{n\ge m}x_{s|2n}-\sum_{n\ge m}x_{t|2m}\|\le\\
\le&\, \sum_{n\ge m}\|x_{s|2n}\|+\|x_{t|2n}\|\le\sum_{n\ge m}32\e_{2n-1}<\delta.
\end{aligned}
$$
On the other hand the lower bound $\|x_{s|2m}-x_{t|2m}\|>\e_{2m}$ supplied by (2) implies
$$\begin{aligned}
\|\sum_{n\in\w}x_{s|2n}-\sum_{n\in\w}x_{t|2n}\|=&\,\|\sum_{n\ge m}x_{s|2n}-\sum_{n\ge m}x_{t|2n}\|\ge\\ \ge&\,\|x_{s|2m}-x_{t|2m}\|-\|\sum_{n>m}(x_{s|2n}-x_{t|2n})\|>\\>&\,\e_{2m}-\sum_{n>m}16\e_{2n-1}>\e_{2m}-32\e_{2m+1}>0
\end{aligned}$$
(the latter two inequalities follow from (i)).
Now we see that $|D_0|=\cc$ and thus $|A\cap (g+B)|\ge|D_0|=\cc$, which implies that
 $\Ind_P(A\cup B)=1$.
\end{proof}

\section{Haar and universally null sets with unit packing index}

In this section in each non-discrete Polish group $G$ we shall construct geometrically large subsets which are small in the sense of measure.

\begin{theorem}\label{t4} Each non-discrete Polish group $G$ contains a closed nowhere dense Haar null subset $C$ with $\Ind_P(C)=1$ and thus $CC^{-1}=G$.
\end{theorem}

\begin{proof} By Theorem~\ref{t3}, the group $G$ contains two closed subsets $A,B\subset G$ with $\ind_P(A)=\ind_P(B)=\cc$ and $\Ind_P(A\cup B)=1$. By Theorem~\ref{megHar}, the sets $A,B$ are Haar null. Then the union $C=A\cup B$ is Haar null, and, being closed in $G$, is nowhere dense.
\end{proof}

Under the set-theoretic assumption $\cov(\mathbb K)=\cc$ in each
non-discrete Polish group $G$ we can construct an universally null
subset $A\subset G$ with $\Ind_P(A)=1$. We recall that a subset
$A$ of a topological space $X$ is called {\em universally null} if
$\mu(A)=0$ for every continuous probability Borel measure $\mu$ on
$X$. A measure $\mu$ on $X$ is {\em continuous} if $\mu(\{x\})=0$
for all $x\in X$. It is clear that each universally null subset of
a non-discrete Polish group in Haar null. 

By $\cov(\mathbb{K})$ we denote the smallest cardinality of a cover of the real line by meager subsets. It is known \cite[19.4]{JW} that the equality $\cov(\mathbb K)=\cc$ is equivalent to MA$_{\mathrm{countable}}$, the Martin Axiom for countable posets.
In particular, $\cov(\mathbb K)=\cc$ holds under MA, the Martin Axiom.

\begin{theorem}\label{unul} Under $\cov(\mathbb{K})= \cc$ each non-discrete
 Polish group $G$ contains a universally null subset $A$ with $\Ind_P(A)=1$.
\end{theorem}

\begin{proof} Fix a countable dense subset $\{x_n\}_{n\in \omega}\subset G$ and
a decreasing sequence $(U_n)_{n\ge 0}$ of open neighborhoods of
the unit $e$ of $G$ with $e=\bigcap_{n\ge 0} U_n$. To each function $f:\w\to\IN$ we can assign a dense $G_\delta$-subset $D_f=\bigcap_{n\in\w}\bigcup_{k\ge n}x_kU_{f(k)}$ of $G$.
Let $\mathbb{N}^\omega=\{f_\alpha :\alpha<\cc\}$ be any enumeration of the
set $\mathbb{N}^\omega$.

The group $G$, being Polish and non-discrete, has size $|G|=\cc$. Let
$G\times G =\{(g_\alpha,g'_\alpha) : \alpha<\cc\}$ be any
enumeration of the product $G\times G$. This enumeration induces an enumeration $G= \{g_\alpha :
\alpha<\cc\}$ of $G$ such that $|\{\alpha: g_\alpha=g\}|=\cc$ for each element
$g\in G$.
Put $G_\alpha=\bigcap_{\beta<\alpha}D_{f_\beta}$ for $\alpha<\cc$.

Observe that the equality $\cov(\mathbb{K})=\cc$ implies that
$|G_\alpha\cap g_\alpha G_\alpha|=\cc$. Otherwise $G$ can
be presented as the union of $<\cc$ many meager subsets:
$$G=(G_\alpha\cap g_\alpha G_\alpha)\cup\bigcup_{\beta<\alpha}(G\setminus D_{f_\beta})\cup\bigcup_{\beta<\alpha}(G\setminus
g_\alpha D_{f_\beta}).
$$

The equality $|G_\alpha\cap g_\alpha G_\alpha|=\cc$, $\alpha<\cc$, allows us to construct inductively a transfinite sequence of points $\{a_\alpha:\alpha<\cc\}\subset G$ such that $$a_\alpha \in G_\alpha\cap g_\alpha G_\alpha\setminus\{a_\beta,g_\beta a_\beta:\beta<\alpha\}$$for all $\alpha<\cc$.

The choice of the enumeration $G=\{g_\alpha:\alpha<\cc\}$ implies that the set $A=\{a_\alpha,g_\alpha a_\alpha:\alpha<\cc\}$ has  $\Ind_P(A)=1$.

It remains to check that $A$ is universally null.
Given a $\sigma$-additive Borel probability
 continuous measure $\mu$ on $G$, for every $n\in\omega$ find a number $f(n)\in
 \mathbb{N}$ such that $\mu(a_nU_{f(n)})<2^{-n}$. It follows that the dense $G_\delta$-subset $D_f$ has measure $\mu(D_f)=0$. Find an ordinal
 $\alpha$ such that $f_\alpha=f$ and 
 observe that $a_\beta,g_\beta a_\beta\in G_{f_\alpha}\subset D_f$ for all $\beta >\alpha$. The inequality $\cov(\mathbb{K})\leq
 \mathrm{non}(\mathbb{L})$ following from the Cichon's diagram  (see \cite{V}) guarantees that  $\mu (\{a_\beta,g_\delta a_\beta:\beta\leq
 \alpha\})=0$ and hence $$\mu(A)\leq \mu (\{a_\beta,g_\beta a_\beta:\beta\leq
 \alpha\})+\mu(D_f).$$
\end{proof}

 \begin{remark} Theorem~\ref{unul} cannot be proved in ZFC because in  Laver's model of ZFC each universally null subset $A$ of a Polish group $G$ has size $|A|<\cc$, which implies that $\ind_P(A)=\Ind_P(A)=\cc$.
\end{remark}

\section{Constructing small subsets with a given (sharp) packing index}

In this section we shall show that Theorem~\ref{megHar} cannot be reversed: nowhere dense
 Haar null sets can have arbitrary packing index.
In fact, we shall construct such sets $A$ with an arbitrary (sharp)
packing index
$$\begin{aligned}
\ind_P^\sharp(A)=&\,\sup\{|S|^+:\mbox{$S\subset G$ and $\{xA\}_{x\in S}$ is disjoint}\},\\
\Ind_P^\sharp(A)=&\,\sup\{|S|^+:\mbox{$S\subset G$ and $\{xA\}_{x\in S}$ is almost disjoint}\}.
\end{aligned}
$$
The formulas
$$
\ind_P(A)=\sup\{\kappa:\kappa<\ind_P^\sharp(A)\}\mbox{ and }
\Ind_P(A)=\sup\{\kappa:\kappa<\Ind_P^\sharp(A)\}
$$show that the sharp packing indices carry more information about a set $A$ comparing to
the usual packing indices.

All possible values of the sharp packing indices of subsets of a
given Abelian group are determined by the following result proved
in \cite{L2}.

\begin{proposition}\label{pL2} Let $G$ be an infinite Abelian group
and $L\subset G$ be a subset with $\Ind_P(L)=1$. For a cardinal
$\kappa\in[2,|G|^+]$ the following conditions are equivalent:
\begin{enumerate}
\item there is a subset $A\subset G$ with $\ind_P^\sharp(A)=\kappa$;
\item there is a subset $A\subset L$ with $\ind^\sharp_P(A)=\Ind^\sharp_P(A)=\kappa$;
\item if $\big|G/[G]_2\big|\le 2$, then $\kappa\ne4$ and if\/ $G=[G]_3$, then $\kappa\ne 3$.
\end{enumerate}
\end{proposition}

Here $[G]_p=\{x\in G:x^p=e\}$ for $p\in\{2,3\}$.

Combining Proposition~\ref{pL2} with Theorems~\ref{t4} and \ref{unul}, we obtain the
main result of this section.

\begin{theorem}\label{t6} Let $G$ be a non-discrete Polish Abelian group and $\kappa$ be
a cardinal such that (i) $2\le\kappa\le|G|^+$, (ii) $k\ne 3$ if $G=[G]_3$, and (iii) $k\ne 4$
if $\big|G/[G]_2\big|\le 2$.
\begin{enumerate}
\item The group $G$ contains a nowhere dense Haar null subset $A$ such that $\ind^\sharp_P(A)=\Ind^\sharp_P(A)=\kappa$;
\item If $\cov(\mathbb{K})=\cc$, then $G$ contains a universally null subset $A$ with
$\ind^\sharp_P(A)=\Ind_P^\sharp(A)=\kappa$.
\end{enumerate}
\end{theorem}

Taking into account that $\ind_P(A)=\sup\{\kappa:\kappa<\ind_P^\sharp(A)\}$, we can apply Theorem~\ref{t6} to deduce the following corollary.

\begin{corollary} Let $G$ be a non-discrete Polish Abelian group and $\kappa$ be a cardinal such that (i) $1\le\kappa\le|G|$, (ii) $k\ne 2$ if $G=[G]_3$, and (iii) $k\ne 3$ if $\big|G/[G]_2\big|\le 2$.
\begin{enumerate}
\item The group $G$ contains a nowhere dense Haar null subset $A$ such that $\ind_P(A)=\Ind_P(A)=\kappa$;
\item If $\cov(\mathbb{K})=\cc$, then $G$ contains a universally null subset $A$ with
$\ind_P(A)=\Ind_P(A)=\kappa$.
\end{enumerate}
\end{corollary}

\section{Acknowledgements}
The authors express their sincere thanks to Heike Mildenberger, Andrzej Ros\l anowski, 
and Lubomyr Zdomskyy for their help in understanding  Shelah's paper \cite{Sh}.

\end{document}